\documentclass[12pt, letterpaper]{article}
\usepackage[utf8]{inputenc}

\title{Inverting the Sweep Map on $(kn,n)$-Dyck Paths: A Simple Algorithm}
\author{Erin Milne \thanks {emilne@champlain.edu}}
\date{Champlain College}

\begin{document}

\begin{titlepage}
\maketitle
\end{titlepage}

\begin{abstract}
We show that our algorithm for inverting the sweep map on $(2n, n)$-Dyck paths works for any $(kn, n)$-Dyck path, where $k$ is an arbitrary positive integer.
\end{abstract}

\section{Introduction to New Work}
In [1], we define the Dyck path and sweep map and show that, if one can recover the level of each south and each west endpoint (denoted by S and W, respectively) from the original Dyck path, then one can invert the sweep map on that Dyck path. We now show that the algorithm we created for finding these levels in the $(2n,n)$-case holds for the $(kn,n)$-case by rewriting our original proofs in terms of an arbitrary positive integer $k$.

\section{Theorems and Proofs for the $(kn, n)$-Case}
We begin by proving a theorem about the levels of south endpoints. This theorem and proof are driectly analogous to those for the $(2n, n)$-case:

\vspace{5mm}

\textbf{Theorem 2.1}: Let $\sigma$ be the output sequence of the sweep map applied from right to left on a $(kn,n)$-Dyck path. If $\sigma_{i}$ = W has level $\tau_{i}$ = $mn$ and $\sigma_{i+1}$ = S, then $\tau_{i+1}$ = $mn$ as well. Moreover, if $\sigma_{i+2}$ = S, then $\tau_{i+2}$ = $mn$ also.

\textbf{Proof}: Let $\sigma_{i}$ = W and $\tau_{i}$ = $mn$. Let $\sigma_{i+1}$ = S and assume $\tau_{i+1}$ = $(m+1)n$. This implies that, as we sweep the Dyck path from right to left, we do not encounter any west ends having level $(m+1)n$ before we encounter $\sigma_{i+1}$ = S having level $(m+1)n$. This means that any west end to the right of $\sigma_{i+1}$ = S has level at most $mn$. Assume there is a west end to the right of $\sigma_{i+1}$ = S with level $mn$. The end immediately to the left of this west end must be a south end with level $(m-k)n$. To the left of this south end there can be at most $k$ consecutive west ends (which takes us back to level $mn$) before there must be another south end, necessarily with level $(m-k)n$ again. But to reach a south end with level $(m+1)n$, we must first reach a west end with level at least $(m+k+1)n$. This west end must ultimately be preceded on the right at some point by west ends with levels $(m+k)n, (m+k-1)n, (m+k-2)n,..., (m+k)n$, contradicting our assumption. Since we are sorting in nondecreasing order, $\sigma_{i+1}$ = S must thus have level $\tau_{i+1}$ = $mn$. This same argument also proves that if $\sigma_{i+2}$ = S, then $\tau_{i+2}$ = $mn$ as well.

Now we reprove our inversion algorithm:

\vspace{5mm}

\textbf{Theorem 2.2}: Inversion Algorithm for the Sweep Map on $(kn,n)$-Dyck Paths:

Input: $\sigma$ (sweep map's output sequence of S's and W's)

Output: A rank sequence $\tau$

Algorithm:

Assign $\tau_{1}=0$

For $i$ = 2 to $kn$:

if $\tau_{i}$ is empty, then:

\hspace{1cm}{ if $\sigma_{i}$ = S, then $\tau_{i}$ = $\tau_{i-1}$ 

\hspace{1cm} else if $\sigma_{i}$ = W, find the number of levels thus far in $\tau$ equal to $\tau_{i-1}$ and then subtract the number of levels thus far in $\tau$ equal to $\tau_{i-1} - kn$ that correspond to S's. If the difference $x$ is positive, assign level $\tau_{i-1} + n$ to the positions 
$\tau$ corresponding to the positions in $\sigma$ of the following $x$ W's (beginning with $\sigma_{i}$ = W). If $x$ is negative, assign level $\tau_{i-1}$ to the positions in $\tau$ corresponding to the positions in $\sigma$ of the following $x$ W's.

else if $\tau_{i}$ is not empty, return.

\textbf{Proof}: We know that the first level is necessarily a 0, and we have already proven that any south end must have the same level as the end preceding it in $\sigma$. Now consider the W's. If a W has level $mn+n$, the east end of the corresponding east step has level $mn$. This east end could be either an S or another W. Thus we find all elements of $\sigma$ that have level $mn$. Now, each of the elements could have a corresponding west end with level $mn+n$ or a corresponding south end with level $mn-kn$. Thus if we subtract the number of S's with level $mn-kn$ from the total number of ends with level $mn$ and get a positive difference $x$, then we know there are $x$ west ends with level $mn+n$. If $x$ is negative, then the next $x$ W's have level $mn$, since $\tau$ is nondecreasing. Since we can fill multiple positions in $\tau$ at once this way, we make sure that our algorithm only looks at positions in $\tau$ that are currently empty as we loop over $i$ to avoid overwriting positions that have already been filled. 

\section{Example}
Now we apply our inversion algorithm to a small example, a $(12, 3)$-Dyck path:

\vspace{5mm}

Let $\sigma$ = SWWSWWWSWWWWWWW. $\tau_{1}$ = 0. $\sigma_{2}$ = W, and so far there is one 0 and no -12's in $\tau$, so $\tau_{2}$ = 3. $\sigma_{3}$ = W, and so far there is one 3 and no -9's in $\tau$, so $\tau_{3}$ = 6. $\sigma_{4}$ = S, so $\tau_{4}$ = 6. $\sigma_{5}$ = W, and so far there are two 6's and and no -6's in $\tau$, so $\tau_{5}$ = 9 and $\tau_{6}$ = 9, since $\sigma_{6}$ = W as well. $\sigma_{7}$ =W, and so far there are two 9's and no -2's in $\tau$, so $\tau_{7}$ = 12 and $\tau_{9}$ = 12, since $\sigma_{9}$ = W as well. $\sigma_{8}$ = S, so $\tau_{8}$ = 12. $\sigma_{10}$ = W, and so far there are three 12's in $\tau$ and one 0 that corresponds to an S, so $\tau_{10}$ = 15 and $\tau_{11}$ = 15, since $\sigma_{11}$ = W as well. $\sigma_{12}$ = W, and so far there are two 15's in $\tau$ and no 3's that correspond to S's, so $\tau_{12}$ = 18 and $\tau_{13}$ = 18, since $\sigma_{13}$ = W as well. $\sigma_{14}$ = W, and so far there are two 18's in $\tau$ and one 6 that corresponds to an S, so $\tau_{14}$ = 21. $\sigma_{15}$ = W, and so far there is one 21 in $\tau$ and no 9's that correspond to S's, so $\tau_{15}$ = 24.\\
Thus $\tau$ = [0 3 6 6 9 9 12 12 12 15 15 18 18 21 24].

\vspace{5mm}

Now that we have $\tau$, we can invert the sweep map. $\tau_{1}$ = 0 corresponds to an S, so we take a north step to level 12. Now we search for the right-most 12 in $\tau$; since we applied the sweep map from right to left, the right-most endpoint with a given level in $\tau$ corresponds to the left-most endpoint with that level on the original path. The right-most 12 corresponds to a W, so we take an east step to level 9. The right-most 9 corresponds to a W, so we take an east step to level 6. The right-most 6 corresponds to an S, so we take a north step to level 18. The right-most 18 corresponds to a W, so we take an east step to level 15. The right-most 15 corresponds to a W, so we take an east step to level 12. The second-to-right-most 12 coresponds to an S, so we take a north step to level 24. The only 24 corresponds to a W, so we take an east step to level 21. The only 21 corresponds to a W, so we take an east step to level 18. The left-most 18 corresponds to a W, so we take an east step to level 15. The left-most 15 corresponds to a W, so we take an east step to level 12. The left-most 12 corresponds to a W, so we take an east step to level 9. The left-most 9 corresponds to a W, so we take an east step to level 6. The left-most 6 corresponds to a W, so we take an east step to level 3. The only 3 corresponds to a W, so we take an east step to level 0 and conclude the path. Thus the original path has step sequence\\
$\sigma^{-1}$ = SWWSWWSWWWWWWWW.

\section{Acknowledgements}
I would like to thank Dr. Gregory Warrington of the University of Vermont for his continued support with this project. I also thank Dr. Melanie Brown and Dr. Scott Stevens of Champlain College, where I am currently employed as an adjunct math instructor, for their moral support.\\

\vspace{5mm}

\begin{center}

\textbf{References}

\end{center}

[1] Erin Milne, Inverting the Sweep Map on (2n,n)-Dyck Paths: A Simple Algorithm, preprint, arXiv:1603.09148

\end{document}